\documentclass[11pt,a4paper]{article}
\usepackage{amsfonts,amsgen,amstext,amsbsy,amsopn,amsfonts,amssymb,amscd}
\usepackage[leqno]{amsmath}
\usepackage[amsmath,amsthm,thmmarks]{ntheorem}
\usepackage{epsf,epsfig}
\usepackage{float}
\usepackage{dsfont}
\usepackage{ebezier,eepic}
\usepackage{color}
\usepackage{tikz}
\usepackage{multirow}
\usepackage{mathrsfs}
\usepackage{graphicx}
\usepackage{subfigure}
\usepackage{epstopdf}
\newtheorem{thm}{Theorem}[section]

\newtheorem{prop}[thm]{Proposition}

\newtheorem{lem}[thm]{Lemma}
\newtheorem{example}[thm]{Example}
\newtheorem{false statement}{False statement}
\newtheorem{cor}[thm]{Corollary}
\newtheorem{fact}[thm]{Fact}

\theoremstyle{definition}
\newtheorem{defn}[thm]{Definition}
\newtheorem{claim}[thm]{Claim}

\newtheorem{conj}[thm]{Conjecture}

\makeatletter \@addtoreset{equation}{section}

\newcommand{\ex}{{\rm ex}}
\def\hh{\mathcal{H}}

\def\hht{\mathcal{T}}

\def\hf{\mathcal{F}}
\def\hg{\mathcal{G}}

\def\ha{\mathcal{A}}
\def\hb{\mathcal{B}}
\def\hs{\mathcal{S}}

\def\ex{\mathbb{E}}
\begin{document}
\title{\bf\Large The maximum sturdiness of intersecting families}
\date{}
\author{Peter Frankl$^1$, Jian Wang$^2$\\[10pt]
$^{1}$R\'{e}nyi Institute, Budapest, Hungary\\[6pt]
$^{2}$Department of Mathematics\\
Taiyuan University of Technology\\
Taiyuan 030024, P. R. China\\[6pt]
E-mail:  $^1$frankl.peter@renyi.hu, $^2$wangjian01@tyut.edu.cn
}
\maketitle

\begin{abstract}
Given a family $\hf\subset 2^{[n]}$ and $1\leq i\neq j\leq n$, we use $\hf(\bar{i},j)$ to denote the family $\{F\setminus \{j\}\colon F\in \hf,\ F\cap \{i,j\}=\{j\}\}$. The sturdiness of $\hf$ is defined as the minimum $|\hf(\bar{i},j)|$ over all $i,j\in [n]$ with $i\neq j$. It has a very natural algebraic definition as well. In the present paper, we consider the maximum sturdiness of $k$-uniform intersecting  families, $k$-uniform $t$-intersecting families and non-uniform $t$-intersecting families.  One of the main results shows that for  $n\geq 36(k+6)$, an intersecting family $\hf\subset \binom{[n]}{k}$ has  sturdiness at most $\binom{n-4}{k-3}$, which is best possible.
\end{abstract}

\section{Introduction}

Let $[n]$ be the standard $n$-set $\{1,2,\ldots,n\}$. We use $2^{[n]}$ to denote the power set of $[n]$ and  $\binom{[n]}{k}$ to denote the collection of all $k$-subsets of $[n]$. For a family $\hf\subset 2^{[n]}$ , let $\hf^c$ be the family of complements, $\{[n]\setminus F\colon F\in \hf\}$. For definiteness let $\hf=\{F_1,F_2,\ldots,F_m\}$.

Let $A(\hf)$ be the {\it incidence matrix} of $\hf$, i.e., the $n$ by $|\hf|$ (i.e., $n$ by $m$) 0-1-matrix with general entry
\[
a_{i\ell} =\left\{\begin{array}{ll}
                 1, & \mbox{ if } i\in F_\ell; \\[5pt]
                 0, & \hbox{ if } i\notin F_\ell.
                \end{array}
              \right.
\]
Note that $A(\hf)+A(\hf^c)$ is the $n$ by $m$ all 1 matrix.

Consider the product $B(\hf) = A(\hf)A(\hf^c)^T$. The general entry of $B=B(\hf)$ is
\begin{align}\label{ineq-1}
b_{ij} =\sum_{1\leq \ell\leq m} a_{i\ell} (1-a_{j\ell}).
\end{align}

For $\hf\subset 2^{[n]}$ and $i,j\in [n]$, let
 \[
 \hf(i) =\{F\setminus \{i\}\colon i\in F\in \hf\},\  \hf(\bar{i}) =\{F\colon i\notin F\in \hf\}.
 \]
Let $\hf(i,\bar{j})$ denote the family $\{F\setminus \{i\}\colon  F\in \hf,\ F\cap \{i,j\}=\{i\}\}$.
From \eqref{ineq-1}, $b_{ii}=0$ is obvious. On the other hand $b_{ij}$ for $i\neq j$ counts the number of $F\in \hf$ with $i\in F$, $j\notin F$. Hence
\begin{align}\label{ineq-2}
b_{ij} = |\hf(i,\bar{j})|.
\end{align}

\begin{defn}
For $\hf\subset 2^{[n]}$  define the {\it sturdiness} $\beta(\hf)$ as $\min_{1\leq i\neq  j \leq n} b_{ij}$.
\end{defn}

It should be clear that for $\hf\subset\tilde{\hf} \subset 2^{[n]}$,
\[
b_{ij}(\hf)\leq b_{ij}(\tilde{\hf}) \mbox{ for all }i,j.
\]
Hence
\begin{align}\label{ineq-3}
\beta(\hf)\leq \beta(\tilde{\hf}).
\end{align}

\begin{prop}[Duality]
\begin{align}\label{ineq-1.1}
\beta(\hf)=\beta(\hf^c).
\end{align}
\end{prop}

\begin{proof}
Let us show that $|\hf(i,\bar{j})|=|\hf^c(\bar{i},j)|$. Let $X=[n]\setminus \{i,j\}$. For any $F\in \hf(i,\bar{j})$,  $F\cup \{i\}\in \hf$. Then $(X\setminus F)\cup \{j\}\in \hf^c$, implying $X\setminus F\in \hf^c(\bar{i},j)$. It is easy to check that $\sigma \colon \hf(i,\bar{j})\mapsto \hf^c(\bar{i},j)$ defined by $\sigma(F)=X\setminus F$ is a bijection. Thus $|\hf(i,\bar{j})|=|\hf^c(\bar{i},j)|$ and \eqref{ineq-1.1} follows.
\end{proof}

For $\hf\subset 2^{[n]}$, $b_{ij}(\hf)=2^{n-2}$ for $i\neq j$, i.e.,
\begin{align}
\beta(2^{[n]}) =2^{n-2}.
\end{align}
Similarly,
\begin{align}
\beta\left(\binom{[n]}{k}\right) =\binom{n-2}{k-1}.
\end{align}

 A family $\hf\subset 2^{[n]}$ is called {\it $t$-intersecting} if $|F\cap F'|\geq t$ for all $F,F'\in \hf$. For $t=1$, we simply say that $\hf$ is intersecting.
 One of the most important results in extremal set theory is the following:

\begin{thm}[The Erd\H{o}s-Ko-Rado Theorem \cite{ekr}]\label{thm-ekr}
 Suppose that $n\geq n_0(k,t)$ and $\hf\subset \binom{[n]}{k}$ is $t$-intersecting. Then
\begin{align}\label{ineq-ekr2}
|\hf| \leq \binom{n-t}{k-t}.
\end{align}
\end{thm}

{\noindent\bf Remark.} For $t=1$ the exact value $n_0(k,t)=(k-t+1)(t+1)$ was proved in \cite{ekr}. For
$t\geq 15$ it is due to \cite{F78}. Finally Wilson \cite{W84} closed the gap $2\leq t\leq 14$ with a proof valid for all $t$.

The corresponding problem for {\it non-uniform} families was solved by Katona.

\begin{thm}[The Katona Theorem \cite{Katona}]\label{thm-katona}
Let $n>t>0$ and suppose that $\hf\subset 2^{[n]}$ is $t$-intersecting. Then (i) or (ii) holds.
\begin{itemize}
  \item[(i)] $n-t$ is even and
  \begin{align}\label{ineq-katona1}
  |\hf| \leq \sum_{i\geq \frac{n+t}{2}}\binom{n}{i}.
  \end{align}
  \item[(ii)] $n-t$ is odd and
  \begin{align}\label{ineq-katona2}
  |\hf| \leq \binom{n-1}{\frac{n+t-1}{2}-1}+\sum_{ i\geq \frac{n+t+1}{2}}\binom{n}{i}.
  \end{align}
\end{itemize}
\end{thm}

Our main results determine the maximum sturdiness of intersecting and $t$-intersecting families.

\begin{thm}\label{thm-main}
Let $n\geq 36(k+6)$ and $\hf\subset \binom{[n]}{k}$ be an intersecting family. Then
\begin{align}
\beta(\hf) \leq \binom{n-4}{k-3}.
\end{align}
\end{thm}

Define the {\it triangle family} as
\[
\hht=\hht(n,k)\colon = \left\{F\in \binom{[n]}{k}\colon |F\cap [3]|\geq 2\right\}
\]
The next claim shows that Theorem \ref{thm-main} is best possible.

\begin{claim}\label{claim-1.6}
For $n\geq 2k$, $\beta(\hht) =\binom{n-4}{k-3}$.
\end{claim}

For $t$-intersecting families, we establish the corresponding best possible result for $n\geq 2(t+3)^2k^2$.

\begin{thm}\label{thm-main-1}
Let $\hf\subset \binom{[n]}{k}$ be a $t$-intersecting family. If $n\geq  2(t+3)^2k^2$ then
\[
\beta(\hf) \leq \binom{n-t-3}{k-t-2}.
\]
\end{thm}

The next two results determine the maximum  sturdiness of a non-uniform $t$-intersecting family.

\begin{thm}\label{thm-2}
If $\hf\subset 2^{[n]}$ is intersecting then $\beta(\hf)\leq 2^{n-3}$.
\end{thm}

\begin{thm}\label{thm-main-3}
Suppose that $\hf\subset 2^{[n]}$ is $t$-intersecting.
\begin{itemize}
  \item[(i)] If $n-t=2s$  then $\beta(\hf) \leq \sum\limits_{
  0\leq j\leq s-1} \binom{n-2}{j}$;
  \item[(ii)]If $n-t=2s+1$ and $n\geq \max\{4(s+2)^2, 36(s+7)\}$, then $\beta(\hf) \leq \binom{n-4}{s-2}+\sum\limits_{
  0\leq j\leq s-1} \binom{n-2}{j}$.
\end{itemize}
\end{thm}

Let us recall a convenient notation. For $\hf\subset \binom{[n]}{k}$ and $A\subset B\subset [n]$, define
\[
\hf(A,B) =\left\{F\setminus B\colon F\in \hf,\ F\cap B=A \right\}.
\]
\section{Shifting and proof of Claim \ref{claim-1.6}}

The most important operation on intersecting families is shifting (cf. \cite{F87} for some basic properties). For $\hf\subset \binom{[n]}{k}$ and $1\leq i<j\leq n$, define the shift
$$S_{ij}(\hf)=\left\{S_{ij}(F)\colon F\in\hf\right\},$$
where
$$S_{ij}(F)=\left\{
                \begin{array}{ll}
                 F':= (F\setminus\{j\})\cup\{i\}, & \mbox{ if } j\in F, i\notin F \text{ and } F'\notin \hf; \\[5pt]
                  F, & \hbox{otherwise.}
                \end{array}
              \right.
$$

A family $\hf\subset \binom{[n]}{k}$ is called {\it shifted}  if $S_{ij}(\hf)=\hf$ holds for all $1\leq i<j\leq n$. It is easy to show  (cf. \cite{F87}) that every $k$-graph can be transformed  into a shifted $k$-graph with the same number of edges by repeated shifting.

Let us define the {\it shifting partial order} $\prec$. For two $k$-sets $A$ and $B$ where $A=\{a_1,\ldots,a_k\}$, $a_1<\ldots<a_k$ and $B=\{b_1,\ldots,b_k\}$, $b_1<\ldots<b_k$ we say that $A$ {\it precedes} $B$ and denote it by $A\prec B$ if $a_i\leq b_i$ for all $1\leq i\leq k$.

\begin{prop}[\cite{F87}]
If $\hf\subset \binom{[n]}{k}$ is shifted, then  $A\prec B$ and $B\in \hf$ imply $A\in \hf$.
\end{prop}

Recall that $\gamma(\hf)=\min_{y\in [n]} |\hf(\bar{y})|$ is called the {\it diversity} of the family. For shifted families both  $\gamma(\hf)$ and $\beta(\hf)$ are easy to compute:
\begin{align}
& \gamma(\hf) = |\hf(\bar{1})|,\label{ineq-2.1}\\[3pt]
& \beta(\hf) =|\hf(\bar{1},n)|. \label{ineq-2.2}
\end{align}

A family $\hf\subset 2^{[n]}$ is called {\it $r$-wise $t$-intersecting} if $|F_1\cap F_2\cap \ldots\cap F_r|\geq t$ for all $F_1,F_2,\ldots,F_r\in \hf$. For $r=2$, we simply say that $\hf$ is $t$-intersecting.

The first important property of shifted families, proved by Erd\H{o}s, Ko and Rado \cite{ekr} is
\begin{flalign}
&\mbox{If $\hf\subset\binom{[n]}{k}$ and $n\geq 2k$ then $\hf(n)\subset \binom{[n-1]}{k-1}$ is intersecting as well.}\\[3pt]
&\mbox{(cf. \cite{F78}) If $\hf$ is shifted and intersecting then $\hf(\bar{1})$ is 2-intersecting.}\label{ineq-new2.4}
\end{flalign}
More generally,
\begin{flalign}\label{ineq-new2.5}
&\mbox{(\cite{F87}) If $\hf$ is shifted and $r$-wise  $t$-intersecting then $\hf(\bar{1})$ is $(t+r-1)$-intersecting.}
\end{flalign}

\begin{cor}\label{cor-key}
If $n\geq (t+r)(k-t-r+2)+2$ and $\hf\subset \binom{[n]}{k}$ is  $r$-wise  $t$-intersecting and shifted, then
\begin{align}
\beta(\hf) \leq \binom{n-t-r-1}{k-t-r}.
\end{align}
\end{cor}
\begin{proof}
Since $\hf$ is shifted, by \eqref{ineq-2.2} we infer that $\beta(\hf) =\hf(\bar{1},n)$. By \eqref{ineq-new2.5}, $\hf(\bar{1})\subset \binom{[2,n]}{k-1}$ is $(t+r-1)$-intersecting. We claim that $\hf(\bar{1},n)$ is also $(t+r-1)$-intersecting. Indeed, if there exist $E_1,E_2\in \hf(\bar{1},n)$ such that $|E_1\cap E_2|\leq t+r-2$ then $E_1\cup \{n\},E_2\cup\{n\}\in \hf(\bar{1})$ and $|E_1\cap E_2|= t+r-2$ follows. Since $|E_1\cup E_2|=2(k-1)-(t+r-2) =2k-t-r<n-2$, there exists $x\in [2,n-1]\setminus E_1\cup E_2$. By shiftedness, $E_2\cup \{x\} \in \hf(\bar{1})$. However, $|(E_1\cup \{n\})\cap (E_2\cup \{x\})|=t+r-2$, contradicting the fact $\hf(\bar{1})$ is $(t+r-1)$-intersecting. Thus $\hf(\bar{1},n)$ is $(t+r-1)$-intersecting. By \eqref{ineq-ekr2} and $n-2\geq (t+r)(k-t-r+2)$, $|\hf(1,\bar{n})|\leq \binom{(n-2)-(t+r-1)}{(k-1)-(t+r-1)}= \binom{n-t-r-1}{k-t-r}$.
\end{proof}

\begin{proof}[Proof of Claim \ref{claim-1.6}]
For $3<i\neq j\leq n$,
\[
b_{ij} = 3\binom{n-5}{k-3}+\binom{n-5}{k-4}= \binom{n-4}{k-3} +2\binom{n-5}{k-3}.
\]
For $1\leq i\leq 3$ and $3<j\leq n$,
\[
b_{ij} =2\binom{n-4}{k-2}+\binom{n-4}{k-3}>\binom{n-4}{k-3}.
\]
For $1\leq j\leq 3$ and $3<i\leq n$,
\[
b_{ij} =\binom{n-4}{k-3}.
\]
For $1\leq i\neq j\leq 3$,
\[
b_{ij} =\binom{n-3}{k-2}>\binom{n-4}{k-3}.
\]
Thus for $n\geq 2k$,
\[
\beta(\hht) =\binom{n-4}{k-3}.
\]
\end{proof}

\section{Proof of Theorem \ref{thm-main}}
\begin{prop}
Let $\hf\subset \binom{[n]}{k}$ be an arbitrary family. Then
\begin{align}\label{ineq-key}
\beta(\hf) \leq \frac{k}{n-1} \gamma(\hf).
\end{align}
\end{prop}

\begin{proof}
Fix $y$ with $\hf(\bar{y}) =\gamma(\hf)$. Then
\[
\sum_{x\in [n]\setminus \{y\}} |\hf(x,\bar{y})|= k|\hf(\bar{y})|.
\]
Hence there exists $x\neq y$ with
\[
|\hf(x,\bar{y})| \leq \frac{k}{n-1}|\hf(\bar{y})| =\frac{k}{n-1} \gamma(\hf).
\]
\end{proof}

\begin{cor}\label{cor-1}
If $\gamma(\hf)<\frac{n-1}{k}\binom{n-4}{k-3} = \frac{n-1}{n-3} \frac{k-2}{k} \binom{n-3}{k-2}$, then $\beta(\hf)<\binom{n-4}{k-3}$.
\end{cor}

For $\{u,v,w\}\subset [n]$, let us introduce the notation:
\[
 \hht_{uvw} =\left\{F\in \binom{[n]}{k}\colon |F\cap \{u,v,w\}|\geq 2\right\}.
\]
We need the following two results.
\begin{thm}[\cite{FW2022}]\label{thm-fw-2022}
Let $n>36k$. Suppose that $\hf\subset \binom{[n]}{k}$ is intersecting. Then
\begin{align}
\gamma(\hf) \leq \binom{n-3}{k-2}.
\end{align}
\end{thm}

\begin{thm}[\cite{FW2024}]\label{thm-fw2024}
Let $\hf\subset \binom{[n]}{k}$ be an intersecting family with $k\geq 3$. Define $\alpha$ by $\gamma(\hf)=\left(1-\alpha\right)\binom{n-3}{k-2}$.  If $0\leq  \alpha<1$ and  $n\geq \frac{36k}{1-\alpha}$, then there exists $\{u,v,w\}\subset [n]$ such that $\hf(\bar{u})=\gamma(\hf)$ and
\begin{align*}
&|\hf\setminus  \hht_{uvw}| \leq  18\alpha \binom{n-33}{k-33}.
\end{align*}
\end{thm}

\begin{proof}[Proof of Theorem \ref{thm-main}]
If $\gamma(\hf)<\frac{k-2}{k} \binom{n-3}{k-2}$, then by Corollary \ref{cor-1} the theorem follows. Thus we may assume that $\gamma(\hf)=(1-\alpha) \binom{n-3}{k-2}$ with $\alpha\leq \frac{2}{k}$.

For $k\geq 3$, $(k+6)(k-2)\geq k^2$. Hence $n\geq 36(k+6)$ and $\alpha\leq \frac{2}{k}$
imply $n\geq \frac{36k}{1-\alpha}$. By Theorem \ref{thm-fw2024} there exists $\{u,v,w\}\subset [n]$ such that $\hf(\bar{u})=\gamma(\hf)$ and
\begin{align}\label{ineq4-1}
&|\hf\setminus \hht_{uvw}| \leq  18\alpha \binom{n-33}{k-33}.
\end{align}

Set $T=\{u,v,w\}$. Let
\begin{align*}
&|\hf(\{u,v\},T)|=f_{uv},\  |\hf(\{u,w\},T)|=f_{uw},\ |\hf(\{v,w\},T)|=f_{vw}
\end{align*}
and
\begin{align*}
&|\hf(\{u\},T)|=g_u,\ |\hf(\{v\},T)|=g_v,\  |\hf(\{w\},T)|=g_w\ \mbox{and } |\hf(\emptyset,T)|=h.
\end{align*}

Then we have
\begin{align*}
&|\hf(\bar{u})| = f_{vw} +g_v+g_w+h =(1-\alpha) \binom{n-3}{k-2},\\[5pt]
&|\hf(\bar{v})| = f_{uw} +g_u+g_w+h \leq   \binom{n-3}{k-2}+g_u+g_w+h,\\[3pt]
&|\hf(\bar{w})| = f_{uv} +g_u+g_v+h \leq  \binom{n-3}{k-2}+g_u+g_v+h.
\end{align*}
It follows that
\begin{align}\label{ineq-2.4}
f_{vw}+f_{uw}+f_{uv}+2(g_u+g_v+g_w)+3h &\leq (3-\alpha) \binom{n-3}{k-2} +2g_u+g_v+g_w+2h \nonumber
\\[3pt]
&\leq (3-\alpha) \binom{n-3}{k-2}+2|\hf\setminus \hf_{uvw}^*|.
\end{align}

For any $x\in [n]\setminus \{u,v,w\}$, let $S=T\cup\{x\}$ and
\[
|\hf(\{x,u,v\},S)|=f_{xuv},\ |\hf(\{x,u\},S)|=g_{xu}, |\hf(\{x\},S)|=h_x.
\]
 Define $f_{xuw},f_{xvw},g_{xv},g_{xw}$ similarly.

Suppose for contradiction that $\beta(\hf)>\binom{n-4}{k-3}$. That is, $|\hf(a,\bar{b})|>\binom{n-4}{k-3}$ for all $a,b\in [n]$. Then for  any $x\in [n]\setminus \{u,v,w\}$,
\begin{align*}
&|\hf(x,\bar{u})| = f_{xvw}+g_{xv}+g_{xw}+h_x >\binom{n-4}{k-3},\\[3pt]
&|\hf(x,\bar{v})| = f_{xuw}+g_{xu}+g_{xw}+h_x >\binom{n-4}{k-3},\\[3pt]
&|\hf(x,\bar{w})| = f_{xuv}+g_{xu}+g_{xv}+h_x >\binom{n-4}{k-3}.
\end{align*}
Adding the above  three inequalities, we get
\[
(f_{xvw}+f_{xuw}+f_{xuv})+2(g_{xv}+g_{xu}+g_{xw})+3h_x>3\binom{n-4}{k-3}.
\]
Summing it over all $x\in [n]\setminus \{u,v,w\}$, we get
\[
(k-2)(f_{vw}+f_{uw}+f_{uv})+2(k-1)(g_u+g_v+g_w)+3k h>3(n-3)\binom{n-4}{k-3}.
\]
It follows that
\[
(f_{vw}+f_{uw}+f_{uv})+2\frac{k-1}{k-2}(g_u+g_v+g_w)+3\frac{k}{k-2} h>3\binom{n-3}{k-2}.
\]
Using \eqref{ineq-2.4}, we get
\[
3\binom{n-3}{k-2}<\frac{2}{k-2}(g_u+g_v+g_w)+\frac{6}{k-2} h+ (3-\alpha) \binom{n-3}{k-2}+2|\hf\setminus \hht_{uvw}|.
\]
By simplifying,
\[
(k-2)\alpha \binom{n-3}{k-2} < 2(g_u+g_v+g_w)+6 h+2(k-2)|\hf\setminus \hht_{uvw}|.
\]
Since $g_u+g_v+g_w+h\leq |\hf\setminus \hht_{uvw}|$,
\[
(k-2)\alpha \binom{n-3}{k-2} < 2(k+1)|\hf\setminus \hht_{uvw}|.
\]
Using \eqref{ineq4-1} we obtain that
\[
(k-2)\alpha \binom{n-3}{k-2} <2(k+1)18\alpha \binom{n-33}{k-33}.
\]
If $\alpha=0$ or $k\leq 32$, then by Theorem \ref{thm-fw2024} $\hf\subset \hht_{uvw}$ and $\beta(\hf)\leq \binom{n-4}{k-3}$, contradicting our assumption. Thus $\alpha \neq 0$ and $k\geq 33$.
Then
\[
(k-2) \binom{n-3}{k-2} <36(k+1) \binom{n-33}{k-33}<36(k+1) \binom{n-33}{k-32}.
\]
Consequently,
\[
\left(\frac{n-3}{k-2}\right)^{30} < 36\frac{k+1}{k-2}\leq 144,
\]
contradicting $n\geq 36(k+6)$. Thus the theorem holds.
\end{proof}

For $n=2k$ Erd\H{o}s, Ko and Rado noted that there are many different intersecting families $\hf\subset \binom{[n]}{k}$ satisfying $|\hf| =\frac{1}{2}\binom{n}{k}=\binom{n-1}{k-1}$. In particular, for $n=6$ there is an intersecting 3-graph $\hg_0\subset\binom{[6]}{3}$ with $|\hg_0|=10$ and $\hg_0$ regular with degree 5. 

Based on $\hg_0$, Huang \cite{huang} proved that Theorem  \ref{thm-fw-2022} does not hold in the range $2k\leq n<(2+\sqrt{3})k$. 

Set $\hf_0=\{F\in \binom{[n]}{k}\colon \mbox{ there exists }G\in \hg_0,\ G\subset F\}$.

\begin{thm}[\cite{huang}]\label{huang}
For $k$ sufficiently large and $3k<n<(2+\sqrt{3})k$,
\[
\gamma(\hf_0) > \binom{n-3}{k-2}.
\]
\end{thm}

For $\hh\subset \binom{[n]}{k}$,  a set $T\subset [n]$ is called a {\it transversal} of $\hh$ if $T\cap H\neq \emptyset$ for all $H\in \hh$.
The {\it transversal number} $\tau(\hh)$ of $\hh$ is defined  as the minimum size of a transversal of  $\hh$.

\begin{prop}\label{prop-1}
 Suppose that  $\hh$ is a regular, intersecting 3-graph on 6 vertices and $|\hh|=10$.
 Then $\tau(\hh)=3$.
\end{prop}

\begin{proof}
By assumption $\hh$ is regular of degree  $(10\times 3)/6=5$. On the other hand, if $P$ is a transversal of $\hh$ of size 2, then all four of its supersets are in $\hh$. The remaining $10-4=6$ edges of $\hh$ all contain (exactly) one of the vertices of $P$. Hence at least one of them has  degree at least 7, a contradiction.
\end{proof}

We need the following well-known fact.

\begin{fact}
Suppose that $\alpha \in (0,1)$ is fixed. Then for fixed $t,\ell$, let $k,n\rightarrow \infty$ with $k/n\rightarrow \alpha$,
\begin{align}\label{ineq-2.5}
\binom{n-t}{k-\ell}/\binom{n}{k} \rightarrow \alpha^\ell(1-\alpha)^{t-\ell}.
\end{align}
\end{fact}

Let us show that for essentially the same range, $\beta(\hf_0)>\binom{n-4}{k-3}$ holds as well.

\begin{prop}
For $k$ sufficiently large and $2k<n<(2+\sqrt{3})k$,
\[
\beta(\hf_0)>\binom{n-4}{k-3}.
\]
\end{prop}
\begin{proof}
Note that $\binom{n-4}{k-3} \rightarrow \alpha^3(1-\alpha) \binom{n}{k}$.
For $6<i\neq j\leq n$,
\begin{align*}
|\hf_0(\bar{i},j)| &= 10\binom{n-8}{k-4}+16\binom{n-8}{k-5}+6\binom{n-8}{k-6}+\binom{n-8}{k-7}\\[3pt]
&\rightarrow \left(10\alpha^4(1-\alpha)^4+16\alpha^5(1-\alpha)^3+6\alpha^6(1-\alpha)^2+\alpha^7(1-\alpha)\right) \binom{n}{k}\\[3pt]
&>\alpha^3(1-\alpha) \binom{n}{k} \mbox{ for } \alpha\in (0.12,0.8).
\end{align*}
For $1\leq i\leq 6$ and $6<j\leq n$,
\begin{align*}
|\hf_0(\bar{i},j)| &= 5\binom{n-7}{k-4}+5\binom{n-7}{k-5}+\binom{n-7}{k-6}\\[3pt]
&\rightarrow \left(5\alpha^4(1-\alpha)^3+5\alpha^5(1-\alpha)^2+\alpha^6(1-\alpha)\right) \binom{n}{k}\\[3pt]
&>\alpha^3(1-\alpha) \binom{n}{k} \mbox{ for } \alpha\in \left(2-\sqrt{3},1\right).
\end{align*}
For $1\leq j\leq 6$ and $6<i\leq n$,
\begin{align*}
|\hf_0(\bar{i},j)|&= 5\binom{n-7}{k-3}+10\binom{n-7}{k-4}+5\binom{n-7}{k-5}+\binom{n-7}{k-6}\\[3pt]
&\rightarrow \left(5\alpha^3(1-\alpha)^4+10\alpha^4(1-\alpha)^3+5\alpha^5(1-\alpha)^2+\alpha^6(1-\alpha)\right) \binom{n}{k}\\[3pt]
&>\alpha^3(1-\alpha) \binom{n}{k} \mbox{ for } \alpha\in (0,1).
\end{align*}
For $1\leq i\neq j\leq 6$,
\begin{align*}
|\hf_0(\bar{i},j)|&= 3\binom{n-6}{k-3}+ 4\binom{n-6}{k-4}+ \binom{n-6}{k-5}\\[3pt]
&\rightarrow \left(3\alpha^3(1-\alpha)^3+4\alpha^4(1-\alpha)^2+\alpha^5(1-\alpha)\right) \binom{n}{k}\\[3pt]
&>\alpha^3(1-\alpha) \binom{n}{k} \mbox{ for } \alpha\in (0,1).
\end{align*}
Thus $\beta(\hf_0)> \binom{n-4}{k-3}$ for $2k<n<(2+\sqrt{3})k$ and $k$ sufficiently large.
\end{proof}

\section{Proof of Theorem \ref{thm-main-1}}

In this section, we determine the maximum  of  the sturdiness of a $t$-intersecting family.

A $t$-intersecting family $\hf\subset \binom{[n]}{k}$ is called {\it saturated} if $\hf$ ceases to be  $t$-intersecting by the addition of any further $k$-sets.
By \eqref{ineq-3} we may always assume that $\hf$ is saturated.

For $\hf\subset \binom{[n]}{k}$, the {\it minimum degree} $\delta(\hf)$ is defined as the minimum of $|\hf(i)|$ over all $i\in [n]$. Let us recall a fundamental result of Huang and Zhao.

\begin{thm}[\cite{HZ}]\label{thm-hz}
For $n>2k$, if $\ha,\hb\subset \binom{[n]}{k}$ are cross-intersecting, then
\[
\delta(\ha)\delta(\hb) \leq \binom{n-2}{k-2}^2.
\]
\end{thm}


Let us use it to estimate sturdiness. For $\hf\subset \binom{[n]}{k}$, a set $T\subset [n]$ is called a   {\it $t$-transversal} of $\hf$ if $|T\cap F|\geq t$ for all $F\in \hf$. Define the {\it $t$-transversal number} $\tau_t(\hf)$ as the minimum size of a $t$-transversal of $\hf$.

\begin{thm}\label{thm-4.2}
Let $\hf\subset \binom{[n]}{k}$ be a  $t$-intersecting family. If $\tau_t(\hf)=t+1$ and $n\geq 2k-t+2$, then
\[
\beta(\hf) \leq \binom{n-t-3}{k-t-2}.
\]
\end{thm}

\begin{proof}
Fix a $t$-transversal $S=\{y_1,y_2,\ldots,y_{t+1}\}$ of size $t+1$. Let
\[
\hg_i=\hf(S\setminus \{y_i\},S) \subset \binom{[n]\setminus S}{k-t}, \ i=1,2,\ldots,t+1.
\]
Since $\hf$ is $t$-intersecting,  $\hg_1,\hg_2,\ldots,\hg_{t+1}$ are pairwise cross-intersecting. By saturatedness, all supersets of $S$ are in $\hf$. Thus,
\[
\hf=\hg_1\cup \hg_2\cup\ldots\cup \hg_{t+1}\cup \{F\in \hf\colon S\subset F\}.
\]
By Theorem \ref{thm-hz}, $\delta(\hg_1)\delta(\hg_2) \leq \binom{n-t-3}{k-t-2}^2$. By symmetry assume $\delta(\hg_1) \leq \binom{n-t-3}{k-t-2}$.
That is, for some $x\in [n]\setminus S$, $|\hg_1(x)|\leq \binom{n-t-3}{k-t-2}$. Then
\[
\beta(\hf) \leq |\hf(x,\overline{y_1})|=|\hg_1(x)| =\delta(\hg_1) \leq \binom{n-t-3}{k-t-2}.
\]
\end{proof}

Recall the definition of the {\it Frankl family}:
\[
\ha_i:=\ha_i(n,k,t)=\left\{F\in \binom{[n]}{k}\colon |F\cap [t+2i]|\geq t+i\right\} \mbox{ for }i=1,2,\ldots, k-t.
\]
One can check that $\beta(\ha_1)=\binom{n-t-3}{k-t-2}$.  Moreover,
\[
\beta(\ha_2)=(t+3)\binom{n-t-5}{k-t-3}+\binom{n-t-5}{k-t-4}=\left((t+3) \frac{n-k-1}{k-t-3}+1\right)\binom{n-t-5}{k-t-4}.
\]
Let  $n-t-4= c(k-t-3)$. Then for $1\leq c<t+2$,
\begin{align*}
\frac{\beta(\ha_2)}{\beta(\ha_1)}&=\left((t+3) \frac{n-k-1}{k-t-3}+1\right) \frac{(k-t-2)(k-t-3)}{(n-t-3)(n-t-4)}\\[3pt]
&=\left((t+3) \frac{n-t-4-(k-t-3)}{k-t-3}+1\right) \frac{(k-t-3)^2}{(n-t-4)^2}\\[3pt]
&=\left((t+3)(c-1)+1\right) \frac{1}{c^2}\\[3pt]
&= \frac{t+3}{c}- \frac{t+2}{c^2}>1.
\end{align*}
Thus for $n< (t+2)(k-t-2)+2$, $\beta(\hf)\leq \binom{n-t-3}{k-t-2}$ does  not necessarily  hold for a $t$-intersecting family $\hf\subset \binom{[n]}{k}$.

Let
\[
\tilde{\ha}_1:=\tilde{\ha}_1(n,k,t)=\left\{F\in \binom{[n]}{k}\colon |F\cap [t+2]|= t+1\right\}.
\]
It is easy to see that  $\beta(\ha)=\beta(\ha_1)$ for all $\tilde{\ha}_1\subset \ha\subset \ha_1$. By Corollary \ref{cor-key} we infer that  $\beta(\ha_1)$ is maximal among shifted $t$-intersecting families for $n\geq (t+2)(k-t)$.

%
%
%
%

For $n$ sufficiently large the maximal diversity is known.

\begin{thm}[\cite{F17}]
Suppose that $\hf\subset \binom{[n]}{k}$ is a $t$-intersecting family, $k>t>0$ and $n\geq 2(t+3)^2k^2$.Then
\[
\gamma(\hf) \leq \binom{n-t-2}{k-t-1}.
\]
\end{thm}


We need the notion of basis for saturated $t$-intersecting families. Let $\hht_t(\hf) $ be the family of all $t$-transversals of $\hf$ of sizes at most $k$.  Define the {\it basis} $\hb=\hb(\hf)$ of $\hf$ as  the collection of containment minimal members  in $\hht_t(\hf)$.

The next lemma is easy to prove.

\begin{lem}[\cite{FW2022-0}]\label{lem4-1}
Suppose that $\hf\subset \binom{[n]}{k}$ is a saturated $t$-intersecting family, $n\geq 2k$. Then (i) and (ii) hold.
\begin{itemize}
  \item[(i)] $\hb$ is a $t$-intersecting antichain,
  \item[(ii)] $\hf=\left\{H\in \binom{[n]}{k}\colon \exists B\in \hb, B\subset H\right\}$.
\end{itemize}
\end{lem}

For any $\ell$ with $\tau_t(\hf)\leq \ell \leq k$, let
\[
\hb^{(\ell)} = \left\{B\in \hb\colon |B|=\ell\right\},\ \hb^{(\leq \ell)} = \bigcup_{i=\tau_t(\hf)}^\ell\hb^{(i)} \mbox{ and }\hb^{(\geq \ell)} = \bigcup_{i=\ell}^k\hb^{(i)}
\]


\begin{lem}[\cite{FW2022-0}]\label{lem4-2}
Suppose that $\hf\subset \binom{[n]}{k}$ is a saturated $t$-intersecting family with $\tau_t(\hf)\geq t+1$ and $\hb=\hb_t(\hf)$. Let $r$ be the smallest integer such that $\tau_t(\hb^{(\leq r)})\geq t+1$. Then
\begin{align}\label{ineq-thb-1}
\sum_{r\leq \ell \leq k} \left(\binom{\ell}{t}\ell k^{\ell-t-1}\right)^{-1}|\hb^{(\ell)}|\leq 1.
\end{align}
\end{lem}

\begin{lem}\label{lem4-3}
Let $\hf\subset \binom{[n]}{k}$ be a saturated $t$-intersecting family with $\tau_t(\hf)= t+2$ and $\hb=\hb_t(\hf)$. Let $y\in B_0\in \hb^{(t+2)}$. If $\tau_t(\hb^{(t+2)})\geq t+1$,
then
\[
\left|\hb^{(t+2)}(\bar{y})\right| \leq 4(t+1)(k-t+2).
\]
\end{lem}

\begin{proof}
We prove the lemma by a branching process. During the proof {\it a sequence} $S=(x_1,x_2,\ldots,x_\ell)$ is an ordered sequence of distinct elements of $X$ and we use $\widehat{S}$ to denote the underlying unordered set $\{x_1,x_2,\ldots,x_\ell\}$. In the first stage,  for each of the $t+1$ $t$-subsets $\{x_1,x_2,\ldots,x_t\}\subset B_0\setminus \{y\}$, define a sequence $(x_1,x_2,\ldots,x_t)$. In the second stage, for each sequence $S_t=(x_1,\ldots,x_t)$ of length $t$, by $\tau_t(\hb^{(t+2)})\geq t+1$ there exists $B_1\in \hb^{(t+2)}$ such that $|\widehat{S_t}\cap B_1|<t$. Then replace $S_t$ by the $|B_1\setminus \widehat{S_t}|$ sequences $(x_1,\ldots,x_t,x_{t+1})$ with $x_{t+1}\in B_1\setminus \widehat{S}_t$. In the third stage, for each sequence $S_{t+1}=(x_1,\ldots,x_t,x_{t+1})$ of length $t+1$, by $\tau_t(\hf)\geq t+2$ there exists $F\in \hf$ such that $|\widehat{S_{t+1}}\cap F|<t$. Then replace $S_{t+1}$ by the $|F\setminus \widehat{S_{t+1}}|$ sequences $(x_1,\ldots,x_t,x_{t+1},x_{t+2})$ with $x_{t+2}\in F\setminus \widehat{S_{t+1}}$. Note that there are $\binom{|B_0\setminus \{y\}|}{t}=\binom{t+1}{t}=t+1$ choices for $(x_1,x_2,\ldots,x_t)$. Since $|B_0\cap B_1|\geq t$ and $\widehat{S_t}\subset B_0$, we have $|\widehat{S_t}\cap B_1|\geq t-2$. Thus there are  $|B_1\setminus \widehat{S_t}|\leq 4$ choices for $x_{t+1}$. Similarly $|B_0\cap F|\geq t$ and $\widehat{S_t}\subset B_0$ imply $|\widehat{S_t}\cap F|\geq t-2$. Then there are  $|F\setminus \widehat{S_{t+1}}|\leq k-t+2$ choices for $x_{t+2}$. Thus after three stages there are at most $4(t+1)(k-t+2)$ sequences of length $t+2$. Let $\hs$ be the collection of all these sequences of length $t+2$.

We are left to show that for any $B\in \hb^{(t+2)}(\bar{y})$, there is some $S\in \hs$ such that $B=\widehat{S}$.  Since $|B\cap (B_0\setminus \{y\})|\geq t$, there is some sequence $S_t=(x_1,x_2,\ldots,x_t)$ with $\widehat{S_t}\subset B$ in the first stage.  In the second stage, since $|B\cap B_1|\geq t$ and $|\widehat{S_t}\cap B_1|<t$, there is some $x_{t+1}\in (B_1\cap B)\setminus \widehat{S_t}$. Then  there is a sequence $S_{t+1}=(x_1,x_2,\ldots,x_t,x_{t+1})$ with $\widehat{S_{t+1}}\subset B$. At the third stage,  since $|B\cap F|\geq t$ and $|\widehat{S_{t+1}}\cap F|<t$, there is some $x_{t+2}\in (F\cap B)\setminus \widehat{S_{t+1}}$. Thus there is a sequence $S_{t+2}=(x_1,x_2,\ldots,x_t,x_{t+1},x_{t+2})\in \hs$ with $\widehat{S_{t+2}}= B$. Therefore,
\[
\left|\hb^{(t+2)}(\bar{y})\right|\leq |\hs| \leq 4(t+1)(k-t+2).
\]
\end{proof}

Let $\binom{[n]}{\leq k}$ denote the collection of all subsets of $[n]$ of size at most $k$. For  $\hg \subset \binom{[n]}{\leq k}$, let
\[
\langle \hg \rangle =\left\{F\in \binom{[n]}{k}\colon \mbox{ there exists }G \in \hg \mbox{ such that }G\subset F\right\}.
\]

\begin{proof}[Proof of Theorem \ref{thm-main-1}]
Let $\hf$ be a saturated $t$-intersecting family and let $\hb=\hb(\hf)$ be its basis.  By Theorem \ref{thm-4.2} we may assume $\tau_t(\hf)\geq t+2$. Let $r$ be the smallest integer such that $\tau_t(\hb^{(\leq r)})\geq t+1$.

\vspace{3pt}
{\bf Case 1. } $r\geq t+3$.
\vspace{3pt}

Then there exists $T\in \binom{[n]}{t}$ such that $T$ is a $t$-transversal of $\hb^{(\leq r-1)}$. By \ref{lem4-1} (i), $\hb$ is $t$-intersecting. Thus $T\subset B$ for all  $B\in \hb^{(\leq r-1)}$.
 Fix a $y\in T$.  By \eqref{ineq-thb-1} and $r\geq t+3$,  we have
\begin{align*}
|\hf(\bar{y})| \leq \sum_{r\leq \ell\leq k} |\hb^{(\ell)}|\binom{n-\ell}{k-\ell}\leq \max_{t+3\leq \ell\leq k} \binom{\ell}{t}\ell k^{\ell-t-1}\binom{n-\ell}{k-\ell}.
\end{align*}
Let $f(n,k,\ell,t)=\binom{\ell}{t}\ell k^{\ell-t-1}\binom{n-\ell}{k-\ell}$. Then $n\geq (t+3)^2k^2$ and $\ell\geq t+3$ imply
\[
\frac{f(n,k,\ell+1,t)}{f(n,k,\ell,t)} = \frac{\ell+1}{\ell+1-t} \frac{\ell+1}{\ell} k\frac{k-\ell}{n-\ell} \leq \frac{(t+4)^2}{4(t+3)} \frac{k^2}{n}\leq 1.
\]
It follows that
\begin{align}\label{ineq-4.3}
|\hf(\bar{y})| \leq  \binom{t+3}{t}(t+3) k^2\binom{n-t-3}{k-t-3}.
\end{align}
Note that
\[
\sum_{x\in [n]\setminus\{y\}} |\hf(x,\bar{y})| =k |\hf(\bar{y})|.
\]
By \eqref{ineq-4.3} there exists $x\in [n]\setminus \{y\}$ such that
\[
|\hf(x, \bar{y})| \leq  \frac{k}{n-1}|\hf(\bar{y})| \leq  \frac{k}{n-1}\binom{t+3}{t}(t+3) k^2\binom{n-t-3}{k-t-3}<\frac{(t+3)^4k^4}{6n(n-k)} \binom{n-t-3}{k-t-2}.
\]
For $n\geq (t+3)^2k^2$, we obtain that
\[
\beta(\hf) \leq |\hf(x, \bar{y})| < \frac{1}{3}\binom{n-t-3}{k-t-2}.
\]

\vspace{3pt}
{\bf Case 2. } $r=t+2$.
\vspace{3pt}

Let $S\in \hb^{(t+2)}$ and fix a $y\in S$. Since $\tau_t(\hb^{(t+2)})\geq t+1$ and $\tau_t(\hf)\geq t+2$, by Lemma \ref{lem4-3} we have
\[
\left|\hb^{(t+2)}(\bar{y})\right| \leq 4(t+1)(k-t+2).
\]
Let $X=\cup \hb^{(t+2)}(\bar{y})$. Then
\[
|X|\leq |\hb^{(t+2)}(\bar{y})|(t+2)\leq 4(t+1)(t+2)(k-t+2)<\frac{n}{2}-1.
\]

By  $n\geq 8(t+1)k^2$ we infer that
\begin{align}\label{ineq-4.6}
|\hb^{(t+2)}(\bar{y})|\binom{n-t-3}{k-t-3} \leq 4(t+1)(k-t+2) \frac{k-t-2}{n-k} \binom{n-t-3}{k-t-2} < \frac{1}{2}\binom{n-t-3}{k-t-2}.
\end{align}
Let $\hf' =\langle \hb^{(\geq t+3)}\rangle$.  By \eqref{ineq-thb-1} and $n\geq tk^2$  we have
\begin{align*}
|\hf'(\bar{y})| \leq \sum_{t+3\leq \ell\leq k} |\hb^{(\ell)}|\binom{n-\ell-1}{k-\ell}&\leq \max_{t+3\leq \ell\leq k} \binom{\ell}{t}\ell k^{\ell-t-1}\binom{n-\ell-1}{k-\ell} \\[3pt]
&\leq  \binom{t+3}{t}(t+3) k^2\binom{n-t-4}{k-t-3}\\[3pt]
&\leq  \frac{(t+3)^4 k^2}{6}\binom{n-t-4}{k-t-3}.
\end{align*}
Note that
\[
\sum_{x\in [n]\setminus (X\cup \{y\})} |\hf'(x, \bar{y})| \leq \sum_{x\in [n]\setminus \{y\}} |\hf'(x, \bar{y})| = k |\hf'(\bar{y})|.
\]
Then there exists $x\in [n]\setminus (X\cup \{y\})$ such that for $n\geq 2(t+3)^2k^2$,
\begin{align}\label{ineq-4.7}
|\hf'(x, \bar{y})| \leq \frac{k}{n-|X|-1}|\hf'(\bar{y})|  \leq \frac{2k}{n}|\hf'(\bar{y})| \leq \frac{(t+3)^4k^4}{3n^2}\binom{n-t-3}{k-t-2}<\frac{1}{2}\binom{n-t-3}{k-t-2}.
\end{align}
Adding \eqref{ineq-4.6} and \eqref{ineq-4.7},  we obtain that
\begin{align*}
\beta(\hf) \leq |\hf(x,\bar{y})|<\binom{n-t-3}{k-t-2}.
\end{align*}
\end{proof}

Recall that $\hf\subset \binom{[n]}{k}$ is called $r$-wise $t$-intersecting if $|F_1\cap F_2\cap \ldots\cap F_r|\geq t$ for all $F_1,F_2,\ldots,F_r\in \hf$. If an $r$-wise $t$-intersecting family $\hf$ is not a star, then $\hf$ is $(t+r-2)$-intersecting. Thus we have the following corollary.

\begin{cor}
Let $\hf\subset \binom{[n]}{k}$ be an  $r$-wise  $t$-intersecting family. If $n\geq 2(t+r+1)^2k^2$ then
\begin{align}\label{ineq-4.8}
\beta(\hf) \leq \binom{n-t-r-1}{k-t-r}.
\end{align}
\end{cor}

Let us note that $\ha_1$ shows that \eqref{ineq-4.8} is best possible.

\section{Proof of Theorems \ref{thm-2} and \ref{thm-main-3}}

\begin{fact}
Let $\hf\subset 2^{[n]}$. Then
\begin{align}
\beta(\hf) \leq \frac{n}{4(n-1)}|\hf|.
\end{align}
\end{fact}

\begin{proof}
For every subset $F\in \hf$ there are $|F|(n-|F|)$ choices $x\in F$, $y\notin F$ such that $F$ contributes 1 to $\hf(x,\bar{y})$. As $a(n-a)\leq \frac{n^2}{4}$,
\begin{align}\label{ineq-4.1}
\sum_{x\in [n], y\in [n]\setminus \{x\}} |\hf(x,\bar{y})| \leq \frac{n^2}{4}|\hf|,
\end{align}
i.e., the average of $|\hf(x,\bar{y})|$ is at most $\frac{n^2|\hf|}{4n(n-1)}=\frac{n}{4(n-1)}|\hf|$.

Note that equality holds iff $n$ is even and $\hf\subset \binom{[n]}{n/2}$. For $n=2\ell+1$ the same proof yields
\begin{align*}
\beta(\hf) \leq \frac{\ell(\ell+1)}{(2\ell+1)2\ell} |\hf| =\frac{\ell+1}{2(2\ell+1)} |\hf|.
\end{align*}
\end{proof}

\begin{proof}[Proof of Theorem \ref{thm-2}]
By \eqref{ineq-3}, we may assume $|\hf|=2^{n-1}$. Hence $\hf$ is a filter, that is, $F\subset G\subset [n]$ and $F\in \hf$ imply $G\in \hf$.

Since for each complementary pair $(H,[n]\setminus H)$ exactly one of them is in $\hf$ we can compute \eqref{ineq-4.1} explicitly
\[
\sum_{F\in \hf} |F|(n-|F|) =\frac{1}{2} \sum_{F\in 2^{[n]}} |F|(n-|F|) =\frac{1}{2} \sum_{x\in [n],y\in [n]\setminus \{x\}} |2^{[n]}(x,\bar{y})|=\frac{1}{2} n(n-1)2^{n-2}.
\]
Thus the average is exactly $2^{n-3}$.
\end{proof}

For $A,B\subset [n]$, define the symmetric difference  $A\Delta B$ as $(A\setminus B)\cup (B\setminus A)$ and define the distance between $A$ and $B$ to be
\[
d(A,B)=|A\Delta B|.
\]
 The Hamming ball of center $C\subset [n]$ and radius $r$ is
 \[
 \hb_r(C) =\left\{B\subset [n] \colon d(B,C) \leq r\right\}.
 \]
 A family $\ha\subset 2^{[n]}$ is called a Hamming ball of center $C\subset [n]$  and radius $r$ iff
 \[
  \hb_r(C)\subseteq  \ha \subseteq   \hb_{r+1}(C).
 \]
 For $\ha,\hb\subset 2^{[n]}$, define
 \[
 d(\ha,\hb) =\min\left\{d(A,B)\colon A\in \ha,\ B\in \hb\right\}.
 \]
 For $\ha\subset 2^{[n]}$, the $d$-neighbourhood of $\ha$ is defined as
 \[
 \Gamma_d \ha =\left\{F\subset [n]\colon d(F, \ha) \leq d\right\}.
 \]

Let us recall Harper's theorem \cite{harper}, see \cite{FF81} for a short proof.

 \begin{thm}[Harper's Theorem \cite{harper}]
 Let $\ha\subset 2^{[n]}$ be a non-empty family. Then there exists a Hamming ball $\ha_0$  such that $|\ha_0|=|\ha|$ and $|\Gamma_d \ha |\geq | \Gamma_d \ha_0|$.
 \end{thm}

 Using Harper's Theorem,  Ahlswede and  Katona \cite{AK} proved the following.

\begin{thm}[\cite{AK}]\label{thm-3.4}
Let $1\leq N\leq 2^n$ and let $\ha,\hb\subset 2^{[n]}$ be two families satisfying $|\ha|=N$, $|A\Delta B|\leq w$ for all $A\in \ha$, $B\in \hb$. Let $\ha_0$ be a Hamming ball of center $\emptyset$ with $|\ha_0|=N$. Then
\[
|\hb| \leq \left|\left\{B\subset [n]\colon |B\Delta A| \leq w \mbox{ for all }A\in \ha_0 \right\}\right|.
\]
\end{thm}

One can derive the following theorem from Theorem \ref{thm-3.4}.

 \begin{thm}[\cite{AK}]\label{thm-3.5}
 Let $\ha,\hb\subset 2^{[n]}$ be families satisfying $|A\Delta B|\leq w$ for all $A\in \ha$, $B\in \hb$.
 \begin{itemize}
  \item[(i)] If $w=2s$  then
  \[
  \min\{|\ha|,|\hb|\} \leq \sum\limits_{0\leq j\leq s} \binom{n}{j}.
  \]
  \item[(ii)]If $w=2s+1$ then
   \[
   \min\{|\ha|,|\hb|\} \leq  \binom{n-1}{s}+\sum\limits_{0\leq j\leq s} \binom{n}{j}.
  \]
\end{itemize}
 \end{thm}

 Note that the special case $\ha=\hb$ of Theorem \ref{thm-3.5} is the classical Kleitman's Diameter Theorem (\cite{kleitman}).

By applying Theorem \ref{thm-3.5} we obtain following result.

\begin{thm}\label{thm-main2}
Suppose that $\hf\subset 2^{[n]}$ is a family with $|F\Delta F'|\leq w$ for all $F,F'\in \hf$.
\begin{itemize}
  \item[(i)] If $w=2s$  then $\beta(\hf) \leq \sum\limits_{
  0\leq j\leq s-1} \binom{n-2}{j}$;
  \item[(ii)]If $w=2s+1$ then $\beta(\hf) \leq \sum\limits_{
  0\leq j\leq s-1} \binom{n-2}{j}+\binom{n-3}{s-1}$.
\end{itemize}
\end{thm}

\begin{proof}
 Let $\ha =\hf(1,\bar{2})$ and $\hb=\hf(\bar{1},2)$. Then $|A\Delta B|\leq w-2$ for all $A\in \ha$, $B\in \hb$. Then by Theorem \ref{thm-3.5}
 we infer that for $w=2s$,
 \[
 \beta(\hf) \leq \min\{|\ha|,|\hb|\}\leq \sum\limits_{
  0\leq j\leq s-1} \binom{n-2}{j}.
 \]
 For $w=2s-1$,
  \[
 \beta(\hf) \leq \min\{|\ha|,|\hb|\}\leq\sum\limits_{
  0\leq j\leq s-1} \binom{n-2}{j}+ \binom{n-3}{s-1}.
 \]
\end{proof}

The Hamming ball  $\hb_s(\emptyset)$ shows that the bound on $\beta(\hf)$ in Theorem \ref{thm-main2} (i) is best possible.
Let $\hf=\hb_s(\emptyset)\cup \hht(n,s+1)$. Then $|F\Delta F'|\leq 2s+1$ for all $F,F'\in \hf$ and
\[
\beta(\hf) = \sum\limits_{0\leq j\leq s-1} \binom{n-2}{j}+\binom{n-4}{s-2}.
\]

\begin{conj}
Suppose that $\hf\subset 2^{[n]}$ is a family with $|F\Delta F'|\leq 2s+1$ for all $F,F'\in \hf$. Then for $n\geq 4(s+1)$,
\[
\beta(\hf) \leq \sum\limits_{0\leq j\leq s-1} \binom{n-2}{j}+\binom{n-4}{s-2}.
  \]
\end{conj}

A family $\hg\subset 2^{[n]}$ is called a $u$-union family if $|G\cup G'|\leq  u$ for $G,G'\in \hg$. It is easy to see that $\hg$ is $u$-union iff $\hg^c$ is $(n-u)$-intersecting. Thus the classical Katona Theorem (Theorem \ref{thm-katona}) determines the maximum size of a $u$-union family as well.

Let us give a relatively general example of a $u$-union family for $u=2s+1$, $s\geq 1$.

\begin{example}\label{example-1}
Fix an intersecting family $\hg\subset \binom{[n]}{s+1}$ and define $\hf=\hg\cup \{F\subset [n]\colon |F|\leq  s\}$. For $n\geq 2s+2$ the Erd\H{o}s-Ko-Rado Theorem (Theorem \ref{thm-ekr}) implies
\begin{align}
|\hf| \leq \binom{n-1}{s} +\sum_{0\leq i\leq  s}\binom{n}{i} =2\sum_{0\leq i\leq s} \binom{n-1}{i},
\end{align}
in line with the Katona Theorem (Theorem \ref{thm-katona}). On the other hand,
\[
\gamma(\hf) =\gamma(\hg) +\sum_{0\leq i\leq s}\binom{n-1}{i} \mbox{ and } \beta(\hf) =\beta(\hg) +\sum_{0\leq i\leq s-1}\binom{n-2}{i}.
\]
\end{example}
Thus,
\begin{align}\label{ineq-5.7}
\beta(\hf) \leq \binom{n-4}{s-2}+\sum_{0\leq i\leq s-1}\binom{n-2}{i} \mbox{ fails unless } \beta(\hg)\leq \binom{n-4}{s-2}.
\end{align}
By Theorem \ref{thm-main}, \eqref{ineq-5.7} holds for Example \ref{example-1} if $n\geq 36(s+7)$ or more generally if the $(2s+1)$-union family contains
no members of size exceeding $s+1$.

For the proof of Theorem \ref{thm-u-uninon}, we need the following two results of the first author.

\begin{thm}[\cite{F16}]\label{thm-frankl1}
Let $n,k,t$ be non-negative integers with $n\geq 2k+t$. Suppose that $\hf\subset \binom{[n]}{k+t}$, $\hg\subset \binom{[n]}{k}$ are cross-intersecting. If $\hf$ is $(t+1)$-intersecting and non-empty then
\[
|\hf|+|\hg| \leq 1+\binom{n}{k} -\binom{n-k-t}{k}.
\]
\end{thm}

\begin{thm}[\cite{F78}]\label{thm-frankl2}
Let $\hf\subset \binom{[n]}{k}$ be a $t$-intersecting family with $n>2k-t$. Then
\[
|\hf| \leq \binom{n}{k-t}.
\]
\end{thm}

\begin{thm}\label{thm-u-uninon}
Suppose that $\hf\subset 2^{[n]}$ is $u$-union.
\begin{itemize}
  \item[(i)] If $u=2s$  then $\beta(\hf) \leq \sum\limits_{
  0\leq j\leq s-1} \binom{n-2}{j}$;
  \item[(ii)]If $u=2s+1$ and $n\geq \max\{4(s+2)^2, 36(s+7)\}$, then $\beta(\hf) \leq \sum\limits_{
  0\leq j\leq s-1} \binom{n-2}{j}+\binom{n-4}{s-2}$.
\end{itemize}
\end{thm}

\begin{proof}
Clearly (i) follows from Theorem \ref{thm-main2} (i). It suffices to prove (ii). Let $\hf\subset 2^{[n]}$ be a $(2s+1)$-union family with $\beta(\hf)$ maximal.

If $\hf$ contains
no members of size exceeding $s+1$, then by Theorem \ref{thm-main} and $n\geq 36(s+7)$ we have $\beta(\hf) \leq \sum\limits_{
  0\leq j\leq s-1} \binom{n-2}{j}+\binom{n-4}{s-2}$. Thus we may assume $\max\{|F|\colon F\in \hf\}>s+1$. By \eqref{ineq-3}, without loss of generality we may also assume that $\hf$ is a down set whence $\hf^{(s+2)}\neq \emptyset$.

Note that the $(2s+1)$-union property implies that $\hf^{(s+2)}$, $\hf^{(s)}$ are cross-intersecting and $\hf^{(s+2)}$ is 3-intersecting. By Theorem \ref{thm-frankl1} and $n\geq 2s+2$ we obtain that
\begin{align}\label{ineq-3.8}
|\hf^{(s+2)}|+|\hf^{(s)}| \leq 1+\binom{n}{s}-\binom{n-s-2}{s} <2(s+2)\binom{n-2}{s-1}.
\end{align}

For $3\leq \ell\leq s$, $\hf^{(s+\ell)}$ is $2\ell-1$-intersecting, by Theorem \ref{thm-frankl2} we have
\[
|\hf^{(s+\ell)}| \leq \binom{n}{(s+\ell) -(2\ell-1)} = \binom{n}{s-\ell+1}.
\]
Thus,
\[
\sum_{3\leq \ell\leq s} |\hf^{(s+\ell)}| \leq \binom{n}{s-2}+\binom{n}{s-3} +\ldots +\binom{n}{0}.
\]
Since $n\geq 3s$ implies
\[
\frac{\binom{n}{s-i}}{\binom{n}{s-i-1}} = \frac{n-s+i+1}{s-i} \geq 2,
\]
by $n\geq 12s$ we have
\begin{align}\label{ineq-3.9}
\sum_{3\leq \ell\leq s} |\hf^{(s+\ell)}| \leq 2\binom{n}{s-2} \leq \left(\frac{n}{n-s}\right)^2\frac{2s}{n-s} \binom{n-2}{s-1}< \frac{1}{4}\binom{n-2}{s-1}.
\end{align}

Set $\hh = \hf^{(s)}\cup \hf^{(s+1)}\cup \hf^{(s+2)}$. We choose distinct vertices $x$ and $y$ uniformly at random. Let $X_i=|\hf^{(s+i)}(x,\bar{y})|$ for $i=0,1,2$. Then
\begin{align*}
\ex (X_i)= |\hf^{(s+i)}| \frac{(s+i)(n-s-i)}{n(n-1)}&<|\hf^{(s+i)}| \frac{(s+2)(n-s-2)}{n(n-1)}\\[3pt]
&<|\hf^{(s+i)}| \frac{s+2}{n}.
\end{align*}
By \eqref{ineq-3.8} and $n\geq 4(s+2)^2$ we have
\[
\ex (X_0+X_2) <  \frac{s+2}{n}\left(|\hf^{(s)}|+|\hf^{(s+2)}|\right) \leq \frac{1}{2} \binom{n-2}{s-1}.
\]
Note that  $\hf^{(s+1)}$ is intersecting and  $\hf^{(s)},\hf^{(s+2)}$ are cross 2-intersecting. If $\hf^{(s+1)}$ is non-trivial intersecting, then by Hilton-Milner Theorem we have
\[
|\hf^{(s+1)}| \leq \binom{n-1}{s} -\binom{n-s-2}{s}+1<(s+1)\binom{n-2}{s-1}.
\]
If $\hf^{(s+1)}$ is a star with center $z$, then $\hf^{(s+1)}(z), \hf^{(s+2)}$ is cross-intersecting. Since $\hf^{(s+2)}$ is non-empty, we infer that
\[
|\hf^{(s+1)}| \leq \binom{n-1}{s} -\binom{n-1-(s+2)}{s-1}<(s+2)\binom{n-2}{s-1}.
\]
In both cases, we have $|\hf^{(s+1)}|<(s+2)\binom{n-2}{s-1}$. By $n\geq 4(s+2)^2$ It follows that
\[
\ex X_1 \leq  \frac{s+2}{n}|\hf^{(s+1)}| \leq  \frac{1}{4} \binom{n-2}{s-1}.
\]
Thus $\ex X<  \frac{3}{4} \binom{n-2}{s-1}$. This implies that there  exist $x,y\in [n]$ such that
\begin{align}\label{ineq-3.10}
\hh(x,\bar{y})< \frac{3}{4} \binom{n-2}{s-1}.
\end{align}
Adding \eqref{ineq-3.9}, \eqref{ineq-3.10} and $|\hf^{(i)}(x,\bar{y})|\leq \binom{n-2}{i-1}$ for $i=1,2,\ldots,s-1$, we conclude that
\[
\beta(\hf) \leq |\hf(x,\bar{y})|\leq \sum_{0\leq i\leq s-1}\binom{n-2}{i}< \binom{n-4}{s-2}+\sum_{0\leq i\leq s-1}\binom{n-2}{i}.
\]
\end{proof}

Recall that if  $\hf$ is $t$-intersecting then $\hf^c$ is $(n-t)$-union. By \eqref{ineq-1.1},  Theorem \ref{thm-main-3} follows from Theorem \ref{thm-u-uninon}.

\section{Concluding Remarks}

In the present paper, we mainly considered the maximum sturdiness of $k$-uniform intersecting  families, $k$-uniform $t$-intersecting families and non-uniform $t$-intersecting families.

Let $\hf\subset 2^{[n]}$ be a $(2s+1)$-union family. In Theorem \ref{thm-u-uninon} (ii), we determine the maximum sturdiness of $\hf$ for $n\geq \max\{4(s+2)^2, 36(s+7)\}$. By Example \ref{example-1} and Theorem \ref{huang}, the same does not hold for the range $2(s+1)<n<(2+\sqrt{3})(s+1)$. 

\begin{conj}
Suppose that $\hf\subset 2^{[n]}$ is $(2s+1)$-union. Then for $n\geq 4(s+1)$,
\[
\beta(\hf) \leq \sum_{
  0\leq j\leq s-1} \binom{n-2}{j}+\binom{n-4}{s-2}.
\]
\end{conj}

\begin{conj}
Suppose that $\hf\subset 2^{[n]}$  is a family with $|F\Delta F'|\leq 2s+1$ for all $F,F'\in \hf$. Then for $n\geq 4(s+1)$,
\[
\beta(\hf) \leq \sum_{
  0\leq j\leq s-1} \binom{n-2}{j}+\binom{n-4}{s-2}.
\]
\end{conj}

A related problem is to consider the maximum sturdiness of an IU-family. A family $\hg\subset 2^{[n]}$ is called an IU-family if $\hg$ and $\hg^c$ are both intersecting. Equivalently, an IU-family $\hg\subset 2^{[n]}$ is both intersecting and 1-union.

\begin{conj}
If $\hg\subset 2^{[n]}$ is an IU-family, then
\begin{align}\label{ineq-4.2}
\beta(\hg)\leq 2^{n-4}.
\end{align}
\end{conj}

For $\hg\subset 2^{[n]}$ and $X\subset [n]$, let $\hg_{\mid X}$ denote the family $\{G\cap X\colon G\in \hg\}$.

\begin{fact}
Let $\hg\subset 2^{[n]}$ be an IU-family. If there exists partition $[n]=X\cup Y$ such that $\hg_{\mid X}$ is intersecting, $\hg_{\mid Y}$ is union, then
$\beta(\hg)\leq 2^{n-4}$.
\end{fact}

\begin{proof}
 By \eqref{ineq-3} we may assume first that all $H\subset [n]$ with $H\cap X$ containing a member of $\hg_{\mid X}$ and $H\cap Y$ contained in a member of $\hg_{\mid Y}$ are in $\hg$.
Then $\hg_{\mid X}$ is an upset (filter) and $\hg_{\mid Y}$ is a down-set (complex). Hence for all $x\in X$, $y\in Y$,
\[
|\hg(\bar{x},y)| \leq \left(\frac{1}{2} 2^{|X|-1}\right)\left(\frac{1}{2} 2^{|Y|-1}\right) = \frac{1}{16} 2^n =2^{n-4}, \mbox{ proving }\eqref{ineq-4.2}.
\]
\end{proof}

\end{document}